\documentclass[a4paper,11pt]{article}
\usepackage{latexsym}
\usepackage{amssymb}
\usepackage{amsfonts}
\usepackage{amsmath}
\usepackage{indentfirst}
\usepackage[dvips]{graphicx}
\usepackage{epsfig}
\newtheorem{prop}{Proposition}[section]
\newtheorem{teor}{Theorem}[section]
\newtheorem{lemma}{Lemma}[section]
\newtheorem{cor}{Corollary}[section]

\newcommand{\ninN}{n\in \mathbf{N}}

\newcommand{\cvd}{\hfill $\blacksquare$\bigskip}

\date{}

\author{Filippo Disanto\thanks{Dipartimento di Scienze Matematiche ed
Informatiche, Pian dei Mantellini, 44, 53100, Siena, Italy \; {\tt
disafili@yahoo.it\quad rinaldi@unisi.it}}\and Luca
Ferrari\thanks{Dipartimento di Sistemi e Informatica, viale
Morgagni 65, 50134 Firenze, Italy {\tt ferrari@dsi.unifi.it\quad
pinzani@dsi.unifi.it}}\and Renzo Pinzani$^\ddag$ \and Simone
Rinaldi$^\dag$}

\title{Catalan numbers and relations\footnote{This
work was very poorly supported by MIUR project: \emph{Aspetti
matematici e applicazioni emergenti degli automi e dei linguaggi
formali}.}} 

\begin{document}

\maketitle

\begin{abstract}
We define the notion of a \emph{Catalan pair} (which is a pair of
binary relations $(S,R)$ satisfying certain axioms) with the aim
of giving a common language to most of the combinatorial
interpretations of Catalan numbers. We show, in particular, that
the second component $R$ uniquely determines the pair, and we give
a characterization of $R$ in terms of forbidden configurations. We
also propose some generalizations of Catalan pairs arising from
some slight modifications of (some of the) axioms.
\end{abstract}

\section{Introduction}

A famous exercise of \cite{St1} proposes to the reader to show
that every item of a long list of combinatorial structures
provides a possible interpretation of the well-known sequence of
Catalan numbers. In addition, since its appearance, many new
combinatorial instances of Catalan numbers (in part due to Stanley
as well \cite{St2}) have been presented by several authors
(\cite{BEM,Cl,MM,MaSh,MaSe}, to cite only a few). What makes
Stanley's exercise even more scary is the request for an explicit
bijection for each couple of structures: even the more skillful
and bold student will eventually give up, frightened by such a
long effort.

The motivation of the present work lies in the attempt of making
the above job as easier as possible. We propose yet another
instance of Catalan numbers, by showing that they count pairs of
binary relations satisfying certain axioms. Of course this is not
the first interpretation of Catalan numbers in terms of binary
relations. For instance, a well-known appearance of Catalan
numbers comes from considering the so-called \emph{similarity
relations}; these have been introduced by Fine \cite{F} and
further studied by several authors \cite{GP,M,Sh}. However, what
we claim to be interesting in our setting is that fairly every
known Catalan structure (or, at least, most of the main ones) can
be obtained by suitably interpreting our relations in the
considered framework. From the point of view of our student, this
approach should result in a quicker way to find bijections:
indeed, it will be enough to guess the correct translation of any
two Catalan structures in terms of our binary relations to get, as
a bonus, the desired bijection. We hope to make this statement
much clearer in section \ref{instances}, where, after the
definition of a \emph{Catalan pair} and the proofs of some of its
properties (pursued on sections \ref{def-prop-enum}), we
explicitly describe some representations of Catalan pairs in terms
of well-known combinatorial objects.

The rest of the paper is devoted to show that Catalan pairs are
indeed a concept that deserves to be better investigated. In
section \ref{ReS} we show that any Catalan pair is uniquely
determined by its second component, and we also provide a
characterization of such a component in terms of forbidden
configurations (which, in our case, are forbidden posets). In
addition, we look at what happens when the second component of a
Catalan pair has some specific properties, namely when it
determines a connected posets or a (possibly distributive)
lattice. We also observe that the first component of a Catalan
pair does not uniquely determine the pair itself, and we give a
description of Catalan pairs having the same first component.
Finally, we propose some generalizations of Catalan pairs: in
section \ref{gener1} we see how to modify the axioms in order to
obtain pairs of relations associated with other important integer
sequences, such as Schr\"oder numbers and central binomial
coefficients; moreover we propose a slight, and very natural,
modification of the crucial axiom in the definition of a Catalan
pair and give an account on what this fact leads to.

Throughout the paper, the reader will find a (not at all
exhaustive) series of open problems. We hope they can serve to
stimulate future research on these topics.

\section{Catalan pairs}\label{def-prop-enum}

%
In what follows, given any set $X$, we denote
$\mathcal{D}=\mathcal{D}(X)$ the \emph{diagonal} of $X$, that is
the relation $\mathcal{D}=\{ (x,x)\; |\; x\in X\}$. Moreover, if
$\theta$ is any binary relation on $X$, we denote by
$\overline{\theta}$ the \emph{symmetrization} of $\theta$, i.e.
the relation $\overline{\theta}=\theta \cup \theta^{-1}$.

\subsection{Basic definitions}

Given a set $X$ of cardinality $n$, let $\mathcal{O}(X)$ be the
set of strict order relations on $X$. By definition, this means
that $\theta \in \mathcal{O}(X)$ when $\theta$ is an irreflexive
and transitive binary relation on $X$. In symbols, this means that
$\theta \cap \mathcal{D}=\emptyset$ and $\theta \circ \theta
\subseteq \theta$.

Now let $(S,R)$ be an ordered pair of binary relations on $X$. We
say that $(S,R)$ is a \emph{Catalan pair} on $X$ when the
following axioms are satisfied:
\begin{itemize}
\item[(i)] $S\in \mathcal{O}(X)$; \hfill (\textbf{ord S})

\item[(ii)] $R\in \mathcal{O}(X)$; \hfill (\textbf{ord R})

\item[(iii)] $\overline{R} \cup \overline{S}=X^2 \setminus
\mathcal{D}$; \hfill (\textbf{tot})

\item[(iv)] $\overline{R} \cap \overline{S}=\emptyset$; \hfill
(\textbf{inters})

\item[(v)] $S\circ R\subseteq R$; \hfill (\textbf{comp})
\end{itemize}

\bigskip

\emph{Remarks.}\begin{enumerate}

\item Observe that, since $S$ and $R$ are both strict order
relations, the two axioms (\textbf{tot}) and (\textbf{inters}) can
be explicitly described by saying that, given $x,y\in X$, with $x
\neq y$, exactly one of the following holds: $xSy$, $xRy$, $ySx$,
$yRx$.

\item Axiom (\textbf{comp}) could be reformulated by using strict
containment, i.e. $S\circ R\subset R$. In fact, it is not
difficult to realize that equality cannot hold since $X$ is
finite. However we prefer to keep our notation, thus allowing to
extend the definition of a Catalan pair to the infinite case.

\item From the above axioms it easily follows that $S\cap
S^{-1}=\emptyset$.

\end{enumerate}

\bigskip

In a Catalan pair $(S,R)$, $S$ (resp. $R$) will be referred to as
the \emph{first} (resp. \emph{second}) \emph{component}. Two
Catalan pairs $(S_1 ,R_1 )$ and $(S_2 ,R_2 )$ on the (not
necessarily distinct) sets $X_1$ and $X_2$, respectively, are said
to be \emph{isomorphic} when there exists a bijection $\xi$ from
$X_1$ to $X_2$ such that $xS_1 y$ if and only if $\xi (x)S_2 \xi
(y)$ and $xR_1 y$ if and only if $\xi (x)R_2 \xi (y)$. As a
consequence of this definition, we say that a Catalan pair has
\emph{size} $n$ when it is defined on a set $X$ of cardinality
$n$. The set of isomorphism classes of Catalan pairs of size $n$
will be denoted $\mathcal{C}(n)$. We will be mainly interested in
the set $\mathcal{C}(n)$, even if, in several specific cases, we
will deal with ``concrete" Catalan pairs. However, in order not to
make our paper dull reading, we will use the term
``Catalan pair" when referring both to a specific Catalan pair and
to an element of $\mathcal{C}(n)$. In the same spirit, to mean
that a Catalan pair has size $n$, we will frequently write
``$(S,R)\in \mathcal{C}(n)$", even if $\mathcal{C}(n)$ is a set of
isomorphism classes. In each situation, the context will clarify
which is the exact meaning of what we have written down.


\bigskip

As an immediate consequence of the definition of a Catalan pair
(specifically, from the fact that all the axioms are universal
propositions), the following property holds.

\begin{prop}\label{substr} Let $(S,R)$ be a Catalan pair on $X$. For any $\widetilde{X}\subseteq
X$, denote by $\widetilde{S}$ and $\widetilde{R}$ the restrictions
of $S$ and $R$ to $\widetilde{X}$, respectively. Then
$(\widetilde{S},\widetilde{R})$ is a Catalan pair on
$\widetilde{X}$.
\end{prop}

\subsection{First properties of Catalan pairs}

In order to get trained with the above definition, we start by
giving some elementary properties of Catalan pairs. All the
properties we will prove will be useful in the rest of the paper.

\begin{prop} Given a Catalan pair $(S,R)$, the following
properties hold:
\begin{enumerate}
\item $S\circ R^{-1}\subseteq R^{-1}$;

\item $R\circ S\subseteq R\cup S$;

\end{enumerate}
\end{prop}

\emph{Proof.}
\begin{enumerate}
\item If $xSyR^{-1}z$, then $xSy$ and $zRy$. Since $x$ and $z$ are
necessarily distinct (this follows from axiom (\textbf{inters})),
it must be either $zRx$, $xRz$, $zSx$ or $xSz$. It is then easy to
check that the three cases $xRz$, $zSx$, $xSz$ cannot hold. For
instance, if $xRz$, then $xRzRy$, whence $xRy$, against
(\textbf{inters}) (since, by hypothesis, $xSy$). Similarly, the
reader can prove that both $zSx$ and $xSz$ lead to a
contradiction. Thus $zRx$, i.e. $xR^{-1}z$.

\item Suppose that $xRySz$. Once again, observe that the elements
$x$ and $z$ are necessarily distinct, thus it must be either
$xRz$, $xSz$, $zRx$ or $zSx$. Similarly as above, it can be shown
that neither $zRx$ nor $zSx$ can hold. For instance, in the first
case, from $zRxRy$ we deduce $zRy$, but we have $ySz$ by
hypothesis. The case $zSx$ can be similarly dealt with.\cvd

\end{enumerate}

\emph{Remark.} As a consequence of this proposition, we have that,
in the definition of a Catalan pair, axiom (\textbf{comp}) can be
replaced by:
\begin{equation}\label{overl}
S\circ \overline{R}\subseteq \overline{R}.
\end{equation}
The above property will be useful in the sequel, when we will
investigate the properties of the relation $R$.

\begin{prop}\label{compstar} Let $(S,R)$ be a pair of binary relations on $X$
satisfying axioms (\textbf{ord S}), (\textbf{ord R}),
(\textbf{tot}) and (\textbf{inters}). Then axiom (\textbf{comp})
is equivalent to:
\begin{displaymath}
\overline{S}\circ R\subseteq R\cup S^{-1}.\qquad \qquad
(\textbf{comp*})
\end{displaymath}
\end{prop}

\emph{Proof.}\quad Assume that axiom (\textbf{comp}) holds and let
$x\overline{S}yRz$. Since $x\overline{S}y$, we have two
possibilities: if $xSy$, then $xSyRz$ and $xRz$. Instead, if
$ySx$, then, being also $yRz$, we get that both the cases $xSz$
and $zRx$ cannot occur. Therefore it must be either $zSx$ or
$xRz$, which means that $(x,z)\in R\cup S^{-1}$.

Conversely, assume that condition (\textbf{comp*}) holds, and
suppose that $xSyRz$.We obviously deduce $x\overline{S}yRz$, and
so we have either $xRz$ or $zSx$. If $zSx$, then $zSxSy$, whence
$zSy$, against the hypothesis $yRz$. Therefore it must be
$xRz$.\cvd

\subsection{Catalan pairs are enumerated by Catalan numbers}

To show that the cardinality of $\mathcal{C}(n)$ is given by the
$n$-th Catalan number $C_n$ we will provide a recursive
decomposition for the structures of $\mathcal{C}(n)$. We recall that the
sequence $C_n$ of Catalan numbers starts $1,1,2,5,14,42,\ldots$
(sequence A000108 in \cite{Sl}) and has generating function $\frac{1-\sqrt{1-4x}}{2x}$.

\bigskip

Given two Catalan pairs, say $(S,R)\in \mathcal{C}(n)$ and
$(S',R')\in \mathcal{C}(m)$, suppose that $S$ and $R$ are defined
on $X=\{ x_1 ,\ldots ,x_n \}$, whereas $S'$ and $R'$ are defined
on $Y=\{ y_1 ,\ldots ,y_m \}$, with $X\cap Y=\emptyset$. We define
the \emph{composition} of $(S,R)$ with $(S',R')$ to be the pair of
relations $(S'',R'')$ on the set $\{ z\} \cup X\cup Y$ of
cardinality $n+m+1$, defined by the following properties:
\begin{itemize}
\item[(i)] $S''$ and $R''$, when restricted to $X$, coincide with
$S$ and $R$, respectively;

\item[(ii)] $S''$ and $R''$, when restricted to $Y$, coincides
with $S'$ and $R'$, respectively;

\item[(iii)] for every $x\in X$ and $y\in Y$, it is $xR''y$;

\item[(iv)] for every $x\in X$, it is $xS''z$;

\item[(v)] for every $y\in Y$, it is $zR''y$;

\item[(vi)] no further relation exists among the elements of $\{
z\} \cup X\cup Y$.
\end{itemize}

\noindent For the composition we will use the standard notation, so that
$(S'',R'')=(S,R)\circ (S',R')$.

\bigskip

\noindent\emph{Remarks.} \begin{enumerate}

\item The above definition of composition can be clearly given in
a more compact form by setting $S''=S\cup S'\cup (X\times \{ z\}
)$ and $R''=R\cup R'\cup ((X\cup \{ z\} )\times Y)$.

\item From the above definition it follows that $S''$ is a strict
order relation on $\{ z\} \cup X\cup Y$ and $z$ is a maximal
element of $S''$. Indeed, if $zS''t$, for some $t$, then
necessarily $t\in Y$ (from (iv)), but from (v) we would also have
$zR''t$, against (vi). Similarly, it can be proved that $R''$ is a
strict order relation on $\{ z\} \cup X\cup Y$ and $z$ is a
minimal element of $R''$.

\end{enumerate}

\begin{prop}\label{uno} Let $\alpha =(S,R)\in \mathcal{C}(n)$ and
$\beta =(S',R')\in \mathcal{C}(m)$ be two Catalan pairs as above.
Then $\alpha \circ \beta =(S'',R'')\in \mathcal{C}(n+m+1)$.
\end{prop}

\emph{Proof.}\quad The fact that $S'',R''\in \mathcal{O}(\{ z\}
\cup X\cup Y)$ is stated in remark 2 above. Moreover, if $t,w\in
\mathcal{O}(\{ z\} \cup X\cup Y)$, with $t\neq w$, then the
following cases are possible:
\begin{itemize}
\item both $t$ and $w$ belong to $X$ or $Y$: in this case $(t,w)$
belongs to exactly one among the relations
$S,S^{-1},R,R^{-1},S',(S')^{-1},R',(R')^{-1}$.

\item $t$ belongs to $X$ and $w$ belongs to $Y$: then $tR''w$, and
no further relation exists between $t$ and $w$; the case $t\in Y$
and $w\in X$ can be treated analogously.

\item $t=z$ and $w\in X$: then the only relation between $t$ and
$w$ is $t(S'')^{-1}w$; and similarly, if $w\in Y$, we have only
$tR''w$.
\end{itemize}

As a consequence, we can conclude that $\overline{R''}\cup
\overline{S''}=(\{ z\} \cup X\cup Y)^2 \setminus \mathcal{D}$ and
$\overline{R''}\cap \overline{S''}=\emptyset$.

Finally, suppose that $t(S''\circ R'')w$. If $t,w$ both belong to
$X$ or else to $Y$, then it is immediate to see that $tR''w$.
Otherwise, suppose that $t$ and $w$ are both different from $z$:
then necessarily $t\in X$ and $w\in Y$, and so $tR''w$. Finally,
the cases $t=z$ and $w=z$ cannot occur, as a consequence of remark
2 above. Thus, we can conclude that, in every case, $tR''w$,
whence $S''\circ R'' \subseteq R''$.\cvd

\begin{lemma}\label{S} Given a Catalan pair $(S,R)$ on $X$, let $x,y$ be
two distinct (if any) maximal elements of $S$. Then there exists
no element $t\in X$ such that $tSx$ and $tSy$.
\end{lemma}

\emph{Proof.}\quad If not, since $x$ and $y$ are maximal for $S$,
then necessarily $x\overline{R}y$. If there were an element $t\in
X$ such that $tSx$ and $tSy$, then, from $tSx\overline{R}y$, we
would get $t\overline{R}y$, against the fact that $tSy$.\cvd



Lemma \ref{S} essentially states that the principal ideals
generated by the maximal elements of $X$ (with respect to $S$) are
mutually disjoint.

\begin{prop}\label{due} Let $\gamma =(S'',R'')$ be a Catalan pair of size
$l\geq 1$. Then there exist unique Catalan pairs $\alpha =(S,R)$
and $\beta =(S',R')$ such that $\gamma =\alpha \circ \beta$.
\end{prop}

\emph{Proof.}\quad Suppose that $\gamma$ is defined on $X_l$ of
cardinality $l$ and let $M(S'')$ be the set of the maximal
elements of $S''$. It is clear that $M(S'')\neq \emptyset$, since
$X_l$ is finite. Define the set $\Phi$ to be the set of all
elements of $M(S'')$ which are minimal with respect to $R''$. We
claim that $|\Phi |=1$. Indeed, since the elements of $M(S'')$ are
an antichain of $S''$, then necessarily they constitute a chain of
$R''$, and so the minimum of such a chain is the only element of
$\Phi$. Set $\Phi =\{ x_0 \}$, we can split $X_l$ into three
subsets, $\{ x_0 \}$, $X$ and $Y$, where $X=\{ x\in X_l \; |\;
xS''x_0 \}$ and $Y=\{ x\in X_l \; |\; x_0 R''x \}$. The reader can
easily check that the above three sets are indeed mutually
disjoint. To prove that their union is the whole $X_l$, let $x\in
X_l$ and suppose that $xS''x_0$ does not hold. Since $x_0$ is
maximal for $S''$, then necessarily $x_0\overline{R''}x$. Suppose,
ab absurdo, that $xR''x_0$. Denoting by $y$ the unique (by the
above lemma) element of $M(S'')$ for which $xS''y$, we would have
$xR''x_0 R''y$, and so $xR''y$, a contradiction. Thus we can
conclude that $x_0 R''x$, as desired. Finally, define
$\alpha=(S,R)$ and $\beta=(S',R')$ as the restrictions of
$(S'',R'')$ to the sets $X$ and $Y$, respectively. The fact that
$\alpha$ and $\beta$ are Catalan pairs follows from proposition
\ref{substr}, whereas the proof that $\alpha \circ \beta =\gamma$
is left to the reader. The uniqueness of the above described
decomposition follows from the fact that $|\Phi |=1$, i.e. there
is only one possibility of choosing $x_0$ so that it satisfies the
definition of composition of Catalan pairs. \cvd

\begin{prop} For any $\ninN$, we have:
\begin{equation}\label{Catalan}
|\mathcal{C}(n+1)|=\sum_{k=0}^{n}|\mathcal{C}(k)|\cdot
|\mathcal{C}(n-k)|.
\end{equation}
Since $|\mathcal{C}(0)|=1$, we therefore have that
$|\mathcal{C}(n)|=C_n$, the $n$-th Catalan number.
\end{prop}

\emph{Proof.}\quad By proposition \ref{due}, giving a Catalan pair
of size $n+1$ is the same as giving two Catalan pairs of sizes $k$
and $n-k$, for a suitable $k$. On the other hand, by proposition
\ref{uno} any two Catalan pairs of sizes $k$ and $n-k$ can be
merged into a Catalan pair of size $n+1$. These arguments
immediately imply formula (\ref{Catalan}).\cvd

\section{Combinatorial interpretations of Catalan pairs}\label{instances}

In this section we wish to convince the reader that fairly every
combinatorial structure counted by Catalan numbers can be
interpreted in terms of Catalan pairs. More precisely, we deem
that any Catalan structure can be described using a suitable
Catalan pair $(S,R)$, where $S$ and $R$ are somehow naturally
defined on the objects of the class. To support this statement, we
will take into consideration here five examples, involving rather
different combinatorial objects, such as matchings, paths,
permutations, trees and partitions. For each of them, we will
provide a combinatorial interpretation in terms of Catalan pairs.

\subsection{Perfect noncrossing matchings and Dyck paths}\label{matchings}

Our first example will be frequently used throughout all the
paper. Given a set $A$ of even cardinality, a \emph{perfect
noncrossing matching} of $A$ is a noncrossing partition of $A$
having all the blocks of cardinality 2. There is an obvious
bijection between perfect noncrossing matchings and well formed
strings of parentheses.

A graphical device to represent a perfect noncrossing matching of
$A$ consists of drawing the elements of $A$ as points on a
straight line and join with an arch each couple of corresponding
points in the matching. Using this representation, we can define
the following relations on the set $X$ of arches of a given
perfect noncrossing matching:

\begin{itemize}
\item for any $x,y\in X$, we say that $xSy$ when $x$ is included
in $y$;

\item for any $x,y\in X$, we say that $xRy$ when $x$ is on the
left of $y$.
\end{itemize}

The reader is invited to check that the above definition yields a
Catalan pair $(S,R)$ on the set $X$.

\bigskip

\emph{Example.}\quad Let $X=\{ a,b,c,d,e,f,g\}$, and let $S$ and
$R$ be defined as follows:

\bigskip

$S=\{ (b,a),(f,e),(f,d),(e,d),(g,d)\}$

\smallskip

$R=\{ (a,c),(a,d),(a,e),(a,f),(a,g),(b,c),(b,d),(b,e),(b,f),(b,g),\\
\phantom{\indent R=\{}(c,d),(c,e),(c,f),(c,g),(e,g),(f,g)\} .$

\bigskip

\begin{figure}[!h]
\begin{center}
\includegraphics[scale=0.5]{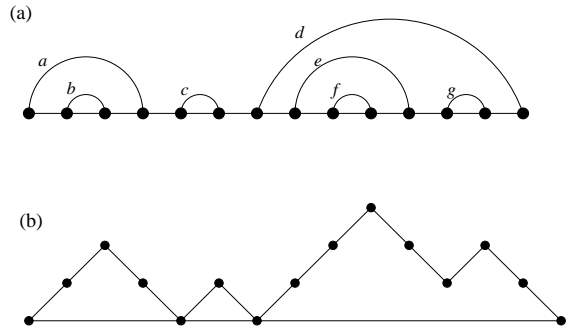}
\end{center}
\caption{The graphical representation of a Catalan pair in terms
of a noncrossing matching, and the associated Dyck
path.}\label{esempio}
\end{figure}

It is easy to check that $(S,R)$ is indeed a Catalan pair on $X$
of size 7, which can be represented as in figure \ref{esempio}(a).

\bigskip

An equivalent way to represent perfect noncrossing matchings is to
use Dyck paths: just interpret the leftmost element of an arch as
an up step and the rightmost one as a down step. For instance, the
matching represented in figure \ref{esempio}(a) corresponds to the
Dyck path depicted in figure \ref{esempio}(b). Coming back to
Catalan pairs, the relations $S$ and $R$ are suitably interpreted
using the notion of tunnel. A \emph{tunnel} in a Dyck path
\cite{E} is a horizontal segment joining the midpoints of an up
step and a down step, remaining below the path and not
intersecting the path anywhere else. Now define $S$ and $R$ on the
set $X$ of the tunnels of a Dyck paths by declaring, for any
$x,y\in X$:
\begin{itemize}
\item $xSy$ when $x$ lies above $y$;

\item $xRy$ when $x$ is completely on the left of $y$.
\end{itemize}

See again figure \ref{esempio} for an example illustrating the
above definition.

\subsection{Pattern avoiding permutations}

Let $n,m$ be two positive integers with $m\leq n$, and let $\pi
=\pi (1)\cdots \pi (n)\in S_n$ and $\nu =\nu (1)\cdots \nu (m)\in
S_m$. We say that $\pi$ {\em contains} the pattern $\nu$ if there
exist indices $i_1 <i_2 <\ldots <i_m$ such that $(\pi ({i_1 }),\pi
({i_2 }),\ldots ,\pi ({i_m }))$ is in the same relative order as
$(\nu (1),\ldots ,\nu(m))$. If $\pi$ does not contain $\nu$, we
say that $\pi$ is {\em $\nu$-avoiding}. See \cite{B} for plenty of
information on pattern avoiding permutations. For instance, if
$\nu =123$, then $\pi =524316$ contains $\nu$, while $\pi =632541$
is $\nu$-avoiding.

We denote by $S_n (\nu )$ the set of $\nu$-avoiding permutations
of $S_n$. It is known that, for each pattern $\nu \in S_3$, $|S_n
(\nu )|=C_n$ (see, for instance, \cite{B}).

\bigskip

It is possible to give a description of the class of 312-avoiding
permutations by means of a very natural set of Catalan pairs. More
precisely, let $[n]=\{ 1,2,\ldots ,n\}$; for every permutation
$\pi \in S_n$, define the following relations $S$ and $R$ on
$[n]$:
\begin{itemize}
\item $iSj$ when $i<j$ and $(j,i)$ is an inversion in $\pi$ (see,
for instance, \cite{B} for the definition of inversion);

\item $iRj$ when $i<j$ and $(i,j)$ is a noninversion in $\pi$.
\end{itemize}

\begin{prop} The permutation $\pi \in S_n$ is 312-avoiding if and
only if $(S,R)$ is a Catalan pair of size $n$.
\end{prop}

\emph{Proof.}\quad The axioms (i) to (iv) in the definition of a
Catalan pair are satisfied by $(S,R)$ for any permutation $\pi$,
as the reader can easily check. Moreover, $\pi$ is 312-avoiding if
and only if, given any three positive integers $i<j<k$, it can
never happen that both $(j,i)$ and $(k,i)$ are inversions and
$(j,k)$ is a noninversion. This happens if and only if $S\circ R$
and $S$ are disjoint. But, from the above definitions of $S$ and
$R$, it must be $S\circ R\subseteq R\cup S$, whence $S\circ
R\subseteq R$.\cvd

The present interpretation in terms of 312-avoiding permutations
can be connected with the previous ones using Dyck paths and
perfect noncrossing matchings, giving rise to a very well-known
bijection, whose origin is very hard to be traced back (see, for
instance, \cite{P}). We leave all the details to the interested
reader.

%

%

\subsection{Plane trees}\label{trees}

By means of the well-known bijection between perfect noncrossing
matchings and plane trees \cite{St1}, the previous example allows
us to give an interpretation of Catalan pairs in terms of plane
trees. The details are left to the reader.

%
%
%
%


\subsection{Noncrossing partitions}

Let $\mathcal{P}_n$ be the set of noncrossing partitions on the
linearly ordered set $X_n =\{ x_1 ,x_2 ,\ldots ,x_n \}$. Each
$p\in \mathcal{P}_n$ determines an equivalence relation $\sim_p$
on $X_n$. Given a generic element $x\in X_n$, we will denote its
equivalence class with $[x]_{\sim_p }$.

Given $x\in X_n$, we set $u(x)=\max_{y<[x]_{\sim_p }}y$. Thus
$u(x)$ is given by the greatest lower bound of the elements in
$[x]_{\sim_p }$ minus 1. Observe that $u(x)$ need not be defined
for all $x$.

Given $p\in \mathcal{P}_n$, define two relations $S$ and $R$ as
follows:

\begin{itemize}

\item $S$ is the transitive closure of the relation $\{ (x,u(x))\;
|\; x\in X_n \}$;

\item $xRy$ when $x<y$ and $(x,y)$ is not in $S$.

\end{itemize}

\begin{figure}[htb]
\begin{center}
\centerline{\hbox{\psfig{figure=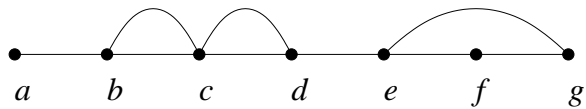,width=3in,clip=}}}
\caption{The noncrossing partition corresponding to the Catalan
pair represented in figure~\ref{esempio}.} \label{noncross}
\end{center}
\end{figure}

Then the pair $(S,R)$ is indeed a Catalan pair on $X_n$, and it
induces an obvious bijection between noncrossing partitions and
plane trees. Figure~\ref{noncross} depicts the noncrossing
partition corresponding to the Catalan pair $(S,R)$ represented in
figure~\ref{esempio}.

\section{Properties of the posets defined by $S$ and $R$}\label{ReS}

In the present section we investigate some features of the posets
associated with the (strict) order relations $R$ and $S$. In the
sequel, a poset will be denoted using square brackets, e.g.
$[X,R]$ and $[X,S]$. An immediate observation which follows
directly from the definition of a Catalan pair is the following,
which we state without proof.

\begin{prop} Given a finite set $X$, consider the graphs $X_1$ and
$X_2$ determined by the Hasse diagrams of the posets $[X,R]$ and
$[X,S]$. Then $X_1$ and $X_2$ are edge-disjoint subgraphs of the
complete graph $K(X)$ on $X$ whose union gives the whole $K(X)$.
\end{prop}


\subsection{The poset defined by $R$}

From the point of view of Catalan pars, it turns out that the
strict order relation $R$ completely defines a Catalan pair. To
prove this, we first need a technical definition which will be
useful again later.

\bigskip

Given a strict order relation $R$ on $X$, define the relation
$\sim_R$ on the set $X$ by declaring $x \sim_R y$ when, for all
$z$, it is $z \overline{R} x$ if and only if $z\overline{R}y$. It
is trivial to show that $\sim_R$ is an equivalence relation. In
what follows, the equivalence classes of $\sim_R$ will be denoted
using square brackets.

\begin{lemma}\label{y}
\begin{enumerate}
\item[(i)] If $x \sim_R y$, then $x\! \! \not{\! R}y$.

\item[(ii)] It is $x \sim_R y$ if and only if, for all $z$, $zRx$
iff $zRy$ and $xRz$ iff $yRz$.

\item[(iii)] If $(S,R)$ is a Catalan pair, then, for all $x,y \in
[z]_{\sim_R}$, it is $xSy$ or $ySx$, i.e. S is a total order on
each equivalence class of $\sim_R$.

\item[(iv)] Suppose $(S,R)$ is a Catalan pair. If $xSy$ and $x\!
\! \not{\! \sim_R}y$, then there exists $a \in X$ such that $a
\overline{R} x$ and $aSy$.

\item[(v)] For all $x,y \in X$, it is $xRy$ iff $[x]_{\sim_R} R
[y]_{\sim_R}$ (where the extension of ${\sim_R}$ to sets has an
obvious meaning).
\end{enumerate}
\end{lemma}

\emph{Proof.} \begin{enumerate}

\item[(i)] Just observe that, if $x\sim_R y$, then $x \overline{R}
y$ would imply $x\overline{R}x$, which is false.

\item[(ii)] Notice that, given that $x\sim_R y$, if $zRx$, then
obviously $z\overline{R}x$, whence $z\overline{R}y$. If we had
$yRz$, then, since $zRx$, it would also be $yRx$, which is
impossible thanks to the preceding statement $(i)$. The fact that
$xRz$ implies $yRz$ can be dealt with analogously.

\item[(iii)] Obvious after $(i)$.

\item[(iv)] From $x\! \! \not{\! \sim_R}y$ it follows, by
definition, that either there exists $a\in X$ such that
$a\overline{R}x$ and $a\! \! \not{\! \overline{R}}y$, or there
exists $b\in X$ such that $b\! \! \not{\! \overline{R}}x$ and
$b\overline{R}y$. The second possibility cannot occur since, if
such an element $b$ existed, then, from the hypothesis $xSy$ and
from (\ref{overl}), we would have $x\overline{R}b$, a
contradiction. Thus an element $a\in X$ with the above listed
properties exists. In particular, since $a\! \! \not{\!
\overline{R}}y$, it must be $a\overline{S}y$. If we had $ySa$,
then, from $xSy$, it would follow $xSa$, a contradiction.
Therefore it must be $aSy$, as desired.

\item[(v)] Suppose that $xRy$. If $a\sim_R x$, applying $(ii)$ it
follows that $aRy$. Now, if it is also $b\sim_R y$, applying
$(ii)$ once more yields $aRb$, which implies the thesis.\cvd
\end{enumerate}

\begin{teor} If $(S_1 ,R),(S_2,R)$ are two Catalan pairs on $X$,
then they are isomorphic.
\end{teor}

\emph{Proof.}\quad From lemma \ref{y}$(iii)$, each equivalence
class of the relation $\sim_R$ is linearly ordered by the order
relations $S_1$ and $S_2$.

Define a function $F$ mapping $X$ into itself such that, if $x \in
X$ and there are exactly $k\geq0$ elements in $[x]_{\sim_R}$ less
than $x$ with respect to the total order $S_1$, then $F(x)$ is
that element in $[x]_{\sim_R}$ having exactly $k$ elements before
it in the total order given by $S_2$.

It is trivial to see that $F$ is a bijection. Since $x \sim_R
F(x)$, using lemma~\ref{y}$(v)$, we get that $xRy$ iff
$F(x)RF(y)$.

To prove that $xS_1y$ implies $F(x)S_2F(y)$ it is convenient to
consider two different cases. First suppose that $x \sim_R y$; in
this case our thesis directly follows from the definition of $F$.
On the other hand, if $x\! \! \not{\! \sim_R}y$, using lemma
\ref{y}$(iv)$, there exists an element $a \in X$ such that $a
\overline{R} x$ and $aS_1y$. Thus, considering the Catalan pair
$(S_2,R)$, it cannot be $F(x) \overline{R} F(y)$, since this would
imply (by lemma \ref{y}$(v)$) that $x \overline{R}y$, against
$xS_1y$. Therefore it must be $F(x) \overline{S_2} F(y)$. More
precisely, we get $F(x)S_2F(y)$, since, from
$F(y)S_2F(x)\overline{R}a$, we would derive $F(y)\overline{R}a$
and so $y\overline{R}a$, which is impossible. With an analogous
argument, we can also prove that $F(x)S_2 F(y)$ implies $xS_1 y$,
which concludes the proof that $F$ is an isomorphism between $(S_1
R)$ and $(S_2 ,R)$.\cvd

%

For the rest of the paper, we set $\mathbf{R}(n)=\{ [X,R]\; |\;
(\exists S)(S,R)\in \mathcal{C}(n)\}$.

\bigskip

\begin{figure}[htb]
\begin{center}
\centerline{\hbox{\psfig{figure=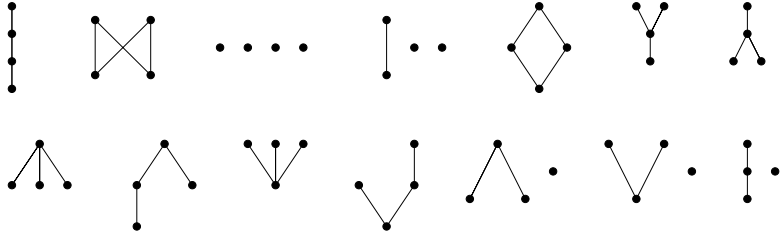,width=4in,clip=}}}
\caption{The 14 posets of $\mathbf{R}(4)$.} \label{erre}
\end{center}
\end{figure}

The posets $[X,R]\in \mathbf{R}(4)$ are those depicted in figure
\ref{erre}.

\begin{figure}[htb]
\begin{center}
\centerline{\hbox{\psfig{figure=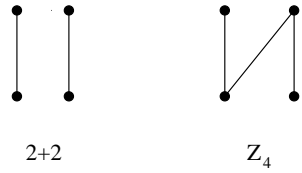,width=1.5in,clip=}}}
\caption{The two posets not belonging to $\mathbf{R}(4)$. }
\label{errenot}
\end{center}
\end{figure}

Among the possible 16 nonisomorphic posets on 4 elements, the two
missing posets are those shown in figure \ref{errenot}. They are
respectively the poset $\mathbf{2}+\mathbf{2}$ (i.e. the direct
sum of two copies of the 2-element chain) and the poset $Z_4$,
called \emph{fence of order 4} (see, for instance,
\cite{C,MZ,St1}).

The rest of this section is devoted to proving that the absence of
the two posets $\mathbf{2}+\mathbf{2}$ and $Z_4$ is not an
accident.

\bigskip

\begin{prop} If $[X,R]\in \mathbf{R}(n)$, then $[X,R]$ does not contain any
subposet isomorphic to $\mathbf{2}+\mathbf{2}$ or $Z_4$.
\end{prop}

\emph{Proof.}\quad Let $(S,R)\in {\mathcal C}(n)$ and suppose, ab
absurdo, that $\mathbf{2}+\mathbf{2}$ is a subposet of $[X,R]$.
Then, denoting with $x,z$ and $y,t$ the minimal and maximal
elements of an occurrence of $\mathbf{2}+\mathbf{2}$ in $[X,R]$,
respectively, and supposing that $xRy$ and $zRt$, we would have,
for instance, $t\overline{S}xRy$. By proposition \ref{compstar},
since $t\! \! \not{\! R}y$, it is $ySt$. However, we also have
$y\overline{S}zRt$ and $y\! \! \not{\! R}t$, whence $tSy$, which
yields a contradiction with the previous derivation.

Similarly, suppose that $Z_4$ is a subposet of $[X,R]$. Then,
supposing that $xRy$, $xRt$ and $zRt$, we have $z\overline{S}xRy$,
whence, by proposition \ref{compstar}, $ySz$. However, it is also
$ySzRt$, which implies $yRt$, and this is false.\cvd

We will now prove that the converse of the above proposition is
also true, thus providing an order-theoretic necessary and
sufficient condition for  a strict order relation $R$ to be the
second component of a Catalan pair.

\begin{prop}\label{sdir} Let $R\in \mathcal{O}(X)$ such that $[X,R]$ does not
contain subposets isomorphic to $\mathbf{2}+\mathbf{2}$ or $Z_4$.
Then $[X,R]\in \mathbf{R}(n)$.
\end{prop}

\emph{Proof.}\quad Given $X=\{ x_1 ,\ldots ,x_n \}$, we define a
binary relation $S=S(R)$ on $X$ by making use of the equivalence
relation $\sim_R$ defined at the beginning of this section. More
precisely:

\begin{enumerate}
\item[-] if $x_i \sim_R x_j$ and $i<j$, set $x_i Sx_j$;

\item[-] if $x\nsim_R y$ and $x\! \! \not{\! \overline{R}}y$, set:
\begin{description}
\item{i)} $xSy$, when there exists $z\in X$ such that
$z\overline{R}x$ and $z\! \! \not{\! \overline{R}}y$;

\item{ii)} $ySx$, when there exists $z\in X$ such that $z\! \!
\not{\! \overline{R}}x$ and $z\overline{R}y$.
\end{description}
\end{enumerate}

We claim that $(S,R)\in \mathcal{C}(n)$.

It is trivial to show that axioms (\textbf{tot}) and
(\textbf{inters}) in the definition of a Catalan pair are
satisfied.

Next we show that axiom (\textbf{comp}) holds. Indeed, suppose
that $xSyRq$ and $x\! \! \! \! \not{\! \!R}q$. From lemma
\ref{y}$(ii)$, it would follow that $x \nsim_R y$. Thus, from
$xSy$ and the definition of $S$, we deduce that there is an
element $z$ such that $z\overline{R}x$ and $z\! \! \! \not{\! \!
\overline{R}}y$. The reader can now check that the four elements
$x,y,q,z$ determine a subposet of $[X,R]$ isomorphic either to
$\mathbf{2}+\mathbf{2}$ or $Z_4$, which is not allowed.

Using an analogous argument, it can be shown that $S \circ R^{-1}
\subseteq R^{-1}$, fact that will be useful below.

Finally, it remains to prove axiom (\textbf{ord S}), i.e. that
$S\in \mathcal{O}(X)$. The fact that $S$ is irreflexive is evident
from its definition. To prove the transitivity of $S$, we first
need to prove that, given $x,y\in X$, the two relations $xSy$ and
$ySx$ cannot hold simultaneously. Indeed, if $x,y\in X$ were such
that $xSy$ and $ySx$, then it could not be $x \sim_R y$ and so, by
definition, there would exist two elements $z,q\in X$ such that
$z\overline{R}x$, $z\! \! \not{\! \overline{R}}y$, $q\! \! \not{\!
\overline{R}}x$ and $q\overline{R}y$. It is not difficult to prove
that the four elements $x,y,z,q$ have to be all distinct (using
the irreflexivity of $R$ and $S$). Now, if we consider the poset
determined by these four elements, in all possible cases a
forbidden poset comes out, and we have reached a contradiction.
Now suppose to have $xSySt$: we want to prove that necessarily
$xSt$. The cases in which we have $x \sim_Ry$ and/or $y \sim_R t$
can be dealt with using the definition of $S$. Moreover, if
$x\nsim_R y$ and $y\nsim_R t$, let $z,q$ such that
$z\overline{R}x$, $z\! \! \! \not{\! \overline{R}}y$, $q\! \! \!
\not{\! \overline{R}}t$ and $q\overline{R}y$. Thanks to the first
part of this proof (namely axiom (\textbf{comp}) and the fact that
$S \circ R^{-1} \subseteq R^{-1}$), from $xSy\overline{R}q$ it
follows that $q\overline{R}x$. On the other hand, if we had
$x\overline{R}t$, since it is $x\! \! \! \not{\! \!
\overline{R}}y$ and $t\! \! \not{\! \overline{R}}y$, it would be
$tSy$ (by the definition of $S$), which is impossible since, by
hypothesis, $ySt$, and we have just shown that the last two
relations lead to a contradiction. Therefore we must have $x\! \!
\not{\! \overline{R}}t$, which, together with $q\! \! \not{\!
\overline{R}}t$ and $q\overline{R}x$, implies that $xSt$, as
desired.\cvd

\begin{figure}[htb]
\begin{center}
\epsfig{file=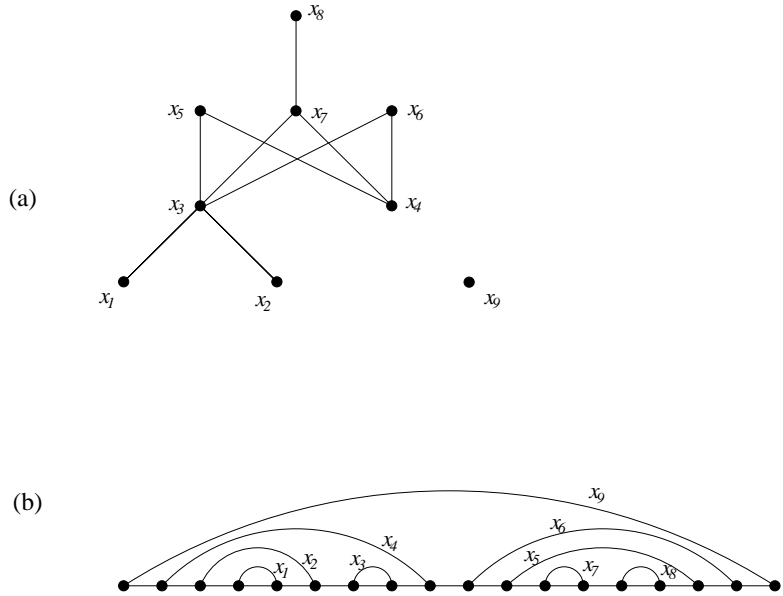, width=4in} \caption{A poset in
$\mathbf{R}(9)$, together with the representation of the
associated Catalan pair as a matching.}\label{s(r)}
\end{center}
\end{figure}

In order to clarify the construction of $S(R)$ given in the proof
of proposition \ref{sdir}, consider the poset $R\in \mathbf{R}(9)$
shown in figure \ref{s(r)}(a). It is $x_1 \sim_R x_2$, hence $x_1
Sx_2$. Similarly we get $x_5 Sx_6$. Moreover, for any fixed
$i=1,\ldots ,8$, we have $x_9 \nsim _R x_i$, and there exists
$x_j$, $j\neq i$, such that $x_i \overline{R}x_j$, so we have $x_i
Sx_9$. Similarly we have $x_2 Sx_4$, $x_3 Sx_4$, $x_7 Sx_5$, $x_7
Sx_6$, $x_8 Sx_5$, $x_8 Sx_6$, and we finally obtain the Catalan
pair $(S,R)$ represented by the matching depicted in figure
\ref{s(r)}(b).

\bigskip

\emph{Remark.}\quad Observe that, as a byproduct of the last
proposition, we have found a presumably new combinatorial
interpretation of Catalan numbers: $C_n$ counts nonisomorphic
posets which are simultaneously $(\mathbf{2}+\mathbf{2})$-free and
$Z_4$-free.

\bigskip

\textbf{Open problem 1.}\quad We have shown that $(S,R)$ is a
Catalan pair if and only if $[X,R]$ does not contain neither
$\mathbf{2}+\mathbf{2}$ nor $Z_4$. The class of
$\mathbf{2}+\mathbf{2}$-free posets have been deeply studied, see
for example \cite{Fis} or the more recent paper \cite{BMCDK}. What
about $Z_4$-free posets?

\bigskip

\textbf{Open problem 2.}\quad Can we define some interesting (and
natural) partial order relation on the set $\mathbf{R}(n)$? Maybe
some of the combinatorial interpretations of Catalan pairs can
help in this task.

\subsection{Imposing some combinatorial conditions on the posets in $\mathbf{R}(n)$}

In this section we impose some conditions on the relation $R$ and
provide the corresponding combinatorial descriptions in terms of
noncrossing matchings and/or Dyck paths and/or 312-avoiding
permutations.

\begin{itemize}

\item[a)] \emph{Connected posets.}\quad First of all notice that a
generic $[X,R]\in \mathbf{R}(n)$ necessarily has at most one
connected component of cardinality greater than one (this follows
at once from the poset avoidance conditions found in the previous
section). It is not difficult to see that, in the interpretation
by means of noncrossing matchings, the fact that $[X,R]$ is
connected means that 1 and $2n$ are not matched. From this
observation it easily follows that $[X,R]$ corresponds to a non
elevated Dyck path and to a 312-avoiding permutation not ending
with 1. This also gives immediately the enumerations of the
Catalan pairs $(S,R)\in \mathcal{C}(n)$ such that $[X,R]$ is
connected. Indeed, elevated Dyck paths of semilength $n$ are known
to be enumerated by the sequence $C_{n-1}$ of shifted Catalan
numbers, whence we get immediately that the number of connected
posets belonging to $\mathbf{R}(n)$ is given by $c_n =C_n
-C_{n-1}$ when $n\geq 2$, whereas $c_0 =c_1 =1$. The resulting
generating function is therefore
\begin{displaymath}
\frac{1-x+2x^2 -(1-x)\sqrt{1-4x}}{2x}.
\end{displaymath}

\begin{figure}[htb]
\begin{center}
\epsfig{file=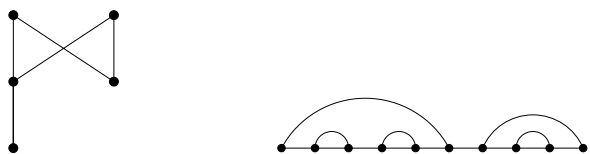, width=3in} \caption{A connected poset
$[X,R]$, and the corresponding perfect matching.}\label{conn}
\end{center}
\end{figure}

\item[b)] \emph{Lattices.}\quad In order to enumerate those posets
of $\mathbf{R}(n)$ which are also lattices, it is convenient to
interpret Catalan pairs as Dyck paths. The following proposition
then holds (where $U$ and $D$ denote up and down steps of Dyck
paths, respectively).

\begin{prop} Let $[X,R]\in \mathbf{R}(n)$ and $P$ be its associated
Dyck path. Then $[X,R]$ is a lattice if and only if $P$ starts and
ends with a peak and does not contain the pattern $DDUU$.
\end{prop}

\emph{Proof.}\quad The fact that $P$ must have a peak both at the
beginning and at the end stems from the fact that a finite lattice
must have a minimum and a maximum. If $P$ contains the pattern
$DDUU$, then denote by $x,y,z,t$ the four tunnels associated with
the four steps of the pattern. It is immediate
to see that $z$ and $t$ are both sups of $x$ and $y$ in $[X,R]$,
which implies that such a poset is not a lattice. Now suppose that
$P$ does not contain the pattern $DDUU$. Given $x,y\in X$
incomparable with respect to $R$, then, in the associated path
$P$, $x$ and $y$ are represented by two tunnels lying one above
the other (say, $x$ above $y$). Consider the down step $D_y$
belonging to $y$. It is obvious that $D_y$ is not isolated, i.e.
it is either followed or preceded by at least another down step.
Now take the first up step coming after $D_y$. Since $P$ avoids
the pattern $DDUU$, such an up step must be followed by a down
step, thus originating a tunnel $z$. It is not difficult to show
that $z$ is the least upper bound of $x$ and $y$. Thus, since any
two elements of $X$ have a least upper bound, we can conclude that
$[X,R]$ is a lattice, as desired.\cvd

\begin{figure}[htb]
\begin{center}
\epsfig{file=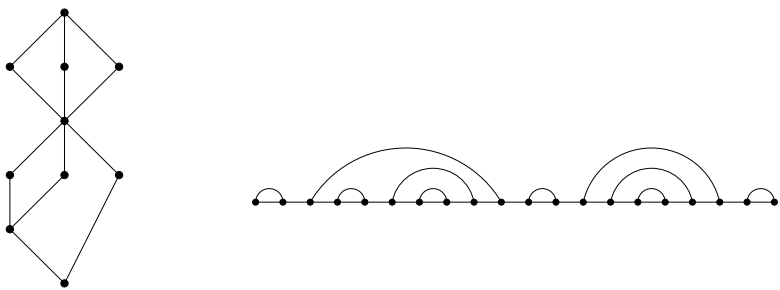, width=4in} \caption{A lattice $[X,R]$, and
the corresponding perfect matching.}\label{lat}
\end{center}
\end{figure}

As a consequence of the last proposition, we are now able to
enumerate Catalan pairs $(S,R)$ such that $[X,R]$ is a lattice.
Indeed, the sequence counting Dyck paths avoiding the pattern
$DDUU$ is A025242 in \cite{Sl} (see also \cite{STT}).

\bigskip

\textbf{Open problem 3.}\quad It seems to be a quite difficult
task to provide a purely order-theoretic characterization of the
lattices $[X,R]$ arising in this way.

\bigskip

\item[c)] \emph{Distributive lattices.}\quad To understand when
$R$ gives rise to a distributive lattice  is undoubtedly a much
easier task. Indeed, in order $[X,R]$ to be a distributive
lattice, it is necessary that it does not contain the two
sublattices $M_3$ and $N_5$~\cite{DP}, shown in figure~\ref{dis}.
This means that, in the associated matching, at most two arches
can be nested and no consecutive sets of nested arches can occur.
Equivalently, the associated Dyck path has height\footnote{The
\emph{height} of a Dyck path is the maximum among the ordinates of
its points.} at most 2, and no consecutive factors of height 2 can
occur. Therefore, an obvious argument shows that the sequence
$d_n$ counting distributive lattices in $\mathbf{R}(n)$ satisfies
the recurrence $d_n =d_{n-1}+d_{n-3}$, with $d_0 =d_1 =d_2 =1$,
having generating function $\frac{1}{1-x-x^3}$, whence $d_n
=\sum_i {n-2i\choose i}$ (sequence A000930 in \cite{Sl}). In this
case, we can also give a structural characterization of
distributive lattices in $\mathbf{R}(n)$: they are all those
expressible as
\begin{displaymath}
\left( \bigoplus_{i=1}^{r}\underline{n_i} \oplus \underline{2}^2
\right)\oplus \underline{n_{r+1}},
\end{displaymath}
where $\oplus$ denotes the \emph{linear} (or \emph{ordinal})
\emph{sum} of posets, $\underline{n}$ is the $n$-element chain and
$\underline{2}^2$ is the Boolean algebra having 4 elements (see
\cite{DP} for basic notions and notations on posets).

\begin{figure}[htb]
\begin{center}
\epsfig{file=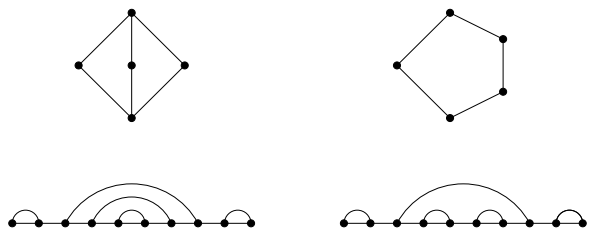, width=3in} \caption{The lattices $M_3$ and $N_5$, and the corresponding perfect matchings.}\label{dis}
\end{center}
\end{figure}

\end{itemize}

\subsection{The poset defined by $S$}

Similarly to what has been done for $R$, we can define the set
$\mathbf{S}(n)=\{ [X,S]\; |\; (\exists R)(S,R)\in
\mathcal{C}(n)\}$. The posets in $\mathbf{S}(n)$ have an
interesting combinatorial characterization, which is described in
the next proposition.

\begin{prop} If $[X,S]\in \mathbf{S}(n)$, then the Hasse diagram
of $[X,S]$ is a forest of rooted trees, where the roots of the
trees are the maximal elements of $S$ and $xSy$ if and only if $y$
is a descendant of $x$ in one of the tree of the forest.
\end{prop}

\emph{Proof.}\quad First observe that, thanks to lemma \ref{S},
the poset $[X,S]$ has $k$ connected components, where $k$ is the
number of its maximal elements. Now take $x,y$ belonging to the
same connected component and suppose that $x\! \! \not{\! S}y$. We
claim that the set of all lower bounds of $\{ x,y\}$ is empty.
Indeed, if we had $z$ such that $zSx$ and $zSy$, then, supposing
(without lose of generality) that $xRy$, it would be $zRy$, a
contradiction. Thus, the Hasse diagram of each connected component
of $[X,S]$ is a direct acyclic graph, that is a tree, rooted at
its maximum element, and this concludes our proof.\cvd

As a consequence of the previous proposition, we have the
following result.

\begin{cor} There is a bijection between $\mathbf{S}(n)$ and the
set of rooted trees with $n+1$ nodes.
\end{cor}

\emph{Proof.}\quad Just add to the Hasse diagram of each element
$[X,S]$ of $\mathbf{S}(n)$ a new root, linking such a root with an
edge to the maximum of each connected component.\cvd

Below the rooted tree on 6 nodes associated with
$(X,S)\in \mathbf{S}(5)$ is shown, where $S=\{ (x_2 ,x_1 ),(x_4
,x_3 ),(x_5 ,x_3 )\}$.

\begin{figure}[htb]
\begin{center}
\centerline{\hbox{\psfig{figure=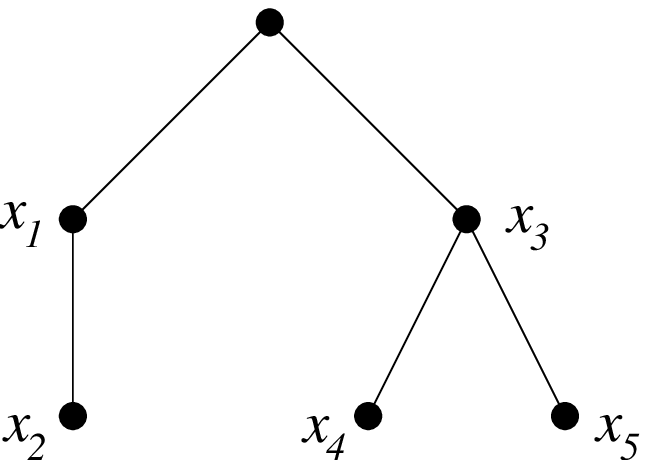,width=1.5in,clip=}}}
\end{center}
\end{figure}

The above corollary implies that $|\mathbf{S}(n)|$ is given by the
number of rooted trees having $n+1$ nodes, which is sequence
A000081 in \cite{Sl}.

\bigskip

Unlike it happens with $R$, the order relation $S$ does not
uniquely determine a Catalan pair. This should be clear by
examining the following two perfect noncrossing
matchings, which are associated with the same $S$, but
determine a different $R$.

\begin{figure}[htb]
\begin{center}
\centerline{\hbox{\psfig{figure=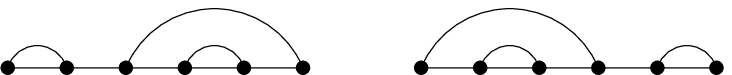,width=3.5in,clip=}}}
\end{center}
\end{figure}

\begin{figure}[htb]
\begin{center}
\centerline{\hbox{\psfig{figure=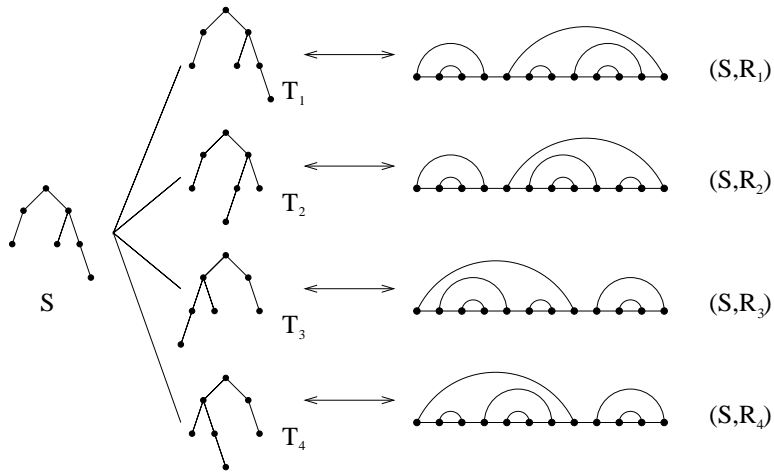,width=4in,clip=}}}
\caption{Representation of the Catalan pairs associated with a
given $S$.} \label{tris}
\end{center}
\end{figure}

This fact is of course an obvious consequence of our last result,
since Catalan pairs are enumerated by Catalan numbers. Recall that
a rooted tree can be seen as a graph-isomorphism class of plane
rooted trees. Since we have shown in section \ref{trees} that
Catalan pairs can be interpreted by using plane rooted trees, it
easily follows that, given $S\in \mathbf{S}(n)$, the set of
Catalan pairs $(S,R)$ can be interpreted as the set of all plane
rooted trees which are isomorphic (as graphs) to the Hasse diagram
of $[X,S]$. Figure \ref{tris} gives an illustration of this
situation, by showing the rooted tree $T$ associated with a given
$S$ and all the plane rooted trees representing the associated
Catalan pairs, together with the alternative representation as
perfect noncrossing matchings.

\section{Generalizations of Catalan pairs}\label{gener1}

In this section we see how a slight modification of the axioms
defining Catalan pairs determines some combinatorial structures
and number sequences, mostly related with permutations. In
particular, we focus our attention on axiom (\textbf{comp}).


We notice that axiom (\textbf{comp}) is the reason since Catalan
pairs can be represented using perfect noncrossing matchings. If
we relax such a condition, we are able to represent some classes
of permutations which, in general, include $312$-avoiding ones.

Consider all pairs of relations $(S,R)$ on a set $X$ satisfying
axioms (\textbf{ord S}), (\textbf{ord R}), (\textbf{tot}) and
(\textbf{inters}). In this situation, we call $(S,R)$ a
\emph{factorial pair} on $X$. The set of all factorial pairs on
$X$ will be denoted $\mathcal{F}(X)$. As we did for Catalan pairs,
we work up to isomorphism, and $\mathcal{F}(n)$ will denote the
isomorphism class of factorial relations on a set $X$ of $n$
elements.

\bigskip

\begin{figure}[htb]
\begin{center}
\centerline{\hbox{\psfig{figure=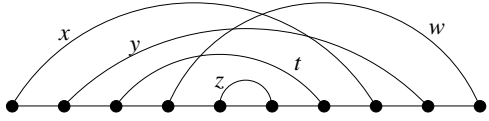,width=2.5in,clip=}}}
\caption{The perfect matching whose associated permutation is
$53124$.}\label{movi+}
\end{center}
\end{figure}

Each pair $(S,R) \in \mathcal{F}(X)$ can be graphically
represented using perfect matchings, extending the encoding given
in section \ref{matchings}. In the matching determined by a
factorial pair, however, two distinct arches can cross, as shown
in figure \ref{movi+}.

The interpretation of the first component of a factorial pair,
$S$, is the same as for Catalan pairs, and corresponds to
\emph{inclusion} of arches. The second component $R$ still
describes the reciprocal position of two arches but, more
generally, we have to consider the reciprocal positions $l(x)$
(left) and $r(x)$ (right) of the two vertices of an arch $x$.
Specifically, we have $xRy$ if and only if $l(x)$ lies on the left
of $l(y)$ and $r(x)$ lies on the left of $r(y)$.

\bigskip

\emph{Example.}\quad Let $(S,R) \in {\cal F}(4)$ represented in
figure \ref{movi+}. Using the notations of figure \ref{movi+}, on
the set of arches $\{ x,y,z,t,w\}$ we have $S=\{
(z,x),(z,y),(z,t),(z,w),(t,x),(t,y)\}$ and $R=\{
(x,y),(x,w),(y,w),(t,w)\}$.

\bigskip

It is clear that, for any set $X$, ${\mathcal C}(X) \subseteq
{\cal F}(X)$. Moreover, using an obvious extension of the
bijection given in section \ref{matchings}, it turns out that
$|{\cal F}(n)|=n!$. More precisely, we have the following
proposition.

\begin{prop}\label{pe} Every factorial pair $(S,R)$ of size $n$ can
be uniquely represented as a permutation $\pi \in S_n$.
\end{prop}

\emph{Proof.}\quad Given $\pi \in S_n$, just define $S$ and $R$ as
in section \ref{matchings}.\cvd



Given a factorial pair $(S,R)$, we call the permutation $\pi$
found in the above proposition its \emph{permutation
representation}. See again figure \ref{movi+} for an example.

\bigskip

Now we come to the main point of the present section, and show how
relaxing axiom (\textbf{comp}) naturally leads to a family of
interesting combinatorial structures which, in some sense,
interpolates between the analogous combinatorial interpretations
of Catalan pairs and factorial pairs.

Denote by $\mathcal{F}_{h,k}(X)$ the class of all pairs of
relations $(S,R)$ on the set $X$ satisfying axioms (\textbf{ord
S}), (\textbf{ord R}), (\textbf{tot}), (\textbf{inters}), and such
that (\textbf{comp}) is replaced by the weaker axiom:

$$ S^h \circ R^k \subseteq R \qquad {\bf (comp \, (h,k) \, )}.$$


The next proposition (whose easy proof is left to the reader)
illustrates how the sets $\mathcal{F}_{h,k}(X)$ are related to
Catalan and factorial pairs.

\begin{prop} \label{}
\begin{enumerate}

\item[(i)] $\mathcal{C}(X) = \mathcal{F}_{1,1}(X)$.

\item[(ii)] For all $h$ and $k$ we have that $\mathcal{F}_{h,k}(X)
\subseteq \mathcal{F}(X)$.

\item[(iii)] If $a \leq b$, then $\mathcal{F}_{a,k}(X) \subseteq
\mathcal{F}_{b,k}(X)$ and $\mathcal{F}_{h,a}(X) \subseteq
\mathcal{F}_{h,b}(X)$ .

\end{enumerate}
\end{prop}

Each element of the family $\{ \mathcal{F}_{h,k}(X) : h,k \geq
1\}$, where $X$ is finite, can be characterized in terms of
permutations avoiding a set of patterns. For example, consider the
two families $\mathcal{F}_{h,1}(X)$ and $\mathcal{F}_{1,k}(X)$.
The following two propositions completely characterize them in
terms of pattern avoiding permutations. The proofs of both
propositions easily follow from the bijection given in proposition
\ref{pe}. In both propositions (as well as in the subsequent
corollary) $X$ denotes a set having $n$ elements.

\begin{prop}\label{av1} The permutation representation of
$\mathcal{F}_{1,k}(X)$ is given by $S_n ((k+2)12\cdots k(k+1))$.
\end{prop}

\begin{prop}\label{av} The permutation representation of $\mathcal{F}_{h,1}(X)$ is
given by $S_n (\pi_2 ,\pi_3 ,\ldots ,\pi_{h+1})$, where $\pi_i \in
S_{h+2}$, for every $2\leq i\leq h+1$, and $\pi_i$ is obtained
from $(h+2)(h+1)\cdots 21$ by moving $i$ to the rightmost
position.
\end{prop}

\begin{cor} \label{omega2} The cardinality of $\mathcal{F}_{2,1}(X)$
is given by the $n$-th Schr\"oder number.
\end{cor}

{\em Proof.}\quad From the previous proposition we get that the
permutation representation of $\mathcal{F}_{2,1}(X)$ is given by
$S_n (4312,4213)$. In \cite{K} it is shown that the above set of
pattern avoiding permutations (or, more precisely, the one
obtained by reversing both patterns) is counted by Schr\"oder
numbers.\cvd

\textbf{Open problem 4.}\quad The enumeration of the sets
$\mathcal{F}_{h,k}(X)$ has to be almost completely carried out,
except for some specific cases. For instance, concerning
$\mathcal{F}_{3,1}(X)$, proposition \ref{av} states that its
permutation representation is given by $S_n (53214,54213,54312)$.
The first terms of its counting sequence are
$1,2,6,24,117,652,3988,\ldots$, which are not in \cite{Sl}.


\section{Other kinds of generalizations}\label{gener1}

Among the possible combinatorial interpretations of Catalan pairs
we have mentioned Dyck paths. In this section we show how some
slight modifications of the axioms for Catalan pairs allow us to
define different pairs of binary relations, which are naturally
interpreted as some well-known families of lattice paths and then
determine well known number sequences. We assume that the reader
is familiar with the most common families of lattice paths, such
as Schr\"{o}der and Grand-Dyck paths.

\bigskip

As usual, we deal with pairs of binary relations $(S,R)$, both
defined on a set $X$ of cardinality $n$ (this will still be
expressed by saying that $(S,R)$ is a pair \emph{of size $n$}) .
The axioms $S$ and $R$ are required to satisfy are the same as the
axioms for Catalan pairs, except for the fact that we do not
impose irreflexivity for $S$. It is immediate to see that all the
remaining axioms are coherent with our new assumption.

\subsection{Unrestricted reflexivity}

Let $\mathcal{U}(n)$ be the set of pairs of binary relations
$(S,R)$ of size $n$, satisfying axioms (\textbf{ord R}),
(\textbf{tot}), (\textbf{inters}), (\textbf{comp}), such that $S$
is a transitive relation and, as it was in the case of Catalan
pairs, $S \cap S^{-1} = \emptyset$. Of course, since we are not
imposing irreflexivity on $S$, given $x\in X$, we may have either
$xSx$ or $x\! \! \not{\! S}x$.

\bigskip

A possible combinatorial interpretation of the elements of
$\mathcal{U}(n)$ can be obtained by means of a slight modification
of the notion of a perfect noncrossing matching. Loosely speaking,
we can introduce two different kinds of arches, namely solid and
dotted arches, imposing that, when $xSx$, the arch corresponding
to $x$ is dotted. These objects will be called \emph{two-coloured
perfect noncrossing matchings} (briefly, \emph{two-coloured
matchings}).

\bigskip

It is evident that, for any Catalan pair $(S,R)\in
\mathcal{C}(n)$, we can define exactly $2^{n}$ different elements
$(S',R')\in \mathcal{U}(n)$ with the property that
\begin{displaymath}
R=R', \quad \mbox{and} \quad S=S'\setminus \mathcal{D}.
\end{displaymath}

\noindent Hence the number of elements of $\mathcal{U}(n)$ is $2^n \, C_n$,
(sequence A052701 in \cite{Sl}).

\bigskip

We obtain some more interesting combinatorial situations by giving
specific axioms for the behavior of $S$ with respect to the
diagonal $\mathcal{D}(X)$.

\subsection{Grand-Dyck paths and central binomial coefficients}

\begin{figure}[htb]
\begin{center}
\epsfig{file=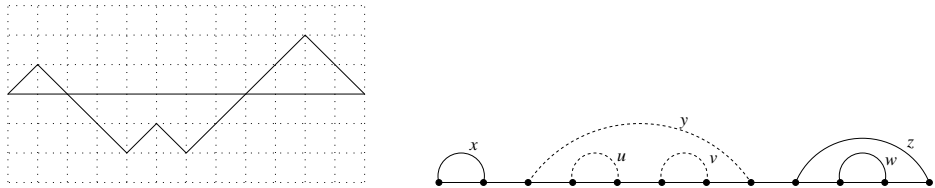, width=4.8in} \caption{A Grand-Dyck path
and its representation as a two-coloured matching.}\label{2eggs}
\end{center}
\end{figure}

Recall that a \emph{Grand-Dyck} path of semi-length $n$ is a
lattice path from $(0,0)$ to $(2n,0)$ using \emph{up} $(1,1)$ and
\emph{down} $(1,-1)$ steps. The number of Grand-Dyck paths of
semi-length $n$ is given by the central binomial coefficient
${2n\choose n}$ \cite{St1}. We can represent a Grand-Dyck path by
using a two-coloured matching, with the convention that for the
parts of the path lying above the $x$-axis we use solid arches,
whereas for the parts of the paths lying below the $x$-axis we use
dotted arches (see figure \ref{2eggs}).

Of course, not every two-coloured matching represents a Grand-Dyck
path. Indeed, we must add the following constraint: if an arch $x$
is contained into an arch $y$, then $x$ and $y$ are either both
solid or both dotted.

In order to give a correct axiomatization of what can be called a
\emph{Grand-Dyck pair}, just add to the axioms for
$\mathcal{U}(n)$ the following one:
\begin{itemize}
  \item if $xSy$, then $xSx$ if and only if $ySy$. \hfill
  (\textbf{choose})
\end{itemize}

Denote by $\mathcal{G}(n)$ the resulting set of pairs of binary
relations, called \emph{Grand-Dyck pairs of size $n$}. It is
evident that, interpreting the relations $S$ and $R$ as in the
case of (one-coloured) matchings, and adding the convention that,
if $xSx$, then $x$ is a dotted arch, we get precisely the set of
two-coloured matchings.

For instance, referring to the example in figure \ref{2eggs}, $R$
and $S$ are as follows:
\begin{displaymath}
\begin{array}{l}
R=\{ (x,y), (x,u), (x,v), (x,z), (y,z), (y,w), (u,v), (u,z),
(u,w),
(v,z), (v,w) \}, \\
S=\{ (u,y), (v,y), (w,z), (u,u), (v,v), (y,y) \}.
\end{array}
\end{displaymath}

Axiom (\textbf{choose}) can be reformulated in a more elegant way.

\begin{prop}\label{alt} Let $\mathcal{D}(S)=\{ (x,x)\in X^2 \; |\; xSx\}$. Then axiom
{\bf (choose)} is equivalent to
\begin{displaymath}
{\cal D}(S) \circ S = S \circ {\cal D}(S).
\end{displaymath}
\end{prop}

\emph{Proof.}\quad Using (\textbf{choose}), it is easy to see that
$x(\mathcal{D}(S)\circ S)y$ if and only if $x\mathcal{D}(S)xSy$ if
and only if $xSy\mathcal{D}(S)y$ if and only if $x(S\circ
\mathcal{D}(S))y$. Conversely, suppose that $xSy$. If $xSx$, then
$x({\cal D}(S)\circ S)y$. However, by hypothesis, this is
equivalent to $x(S\circ {\cal D}(S))y$, whence $ySy$.\cvd

\subsection{Schr\"oder paths and Schr\"oder numbers}

Recall that a \emph{Schr\"oder path} of semi-length $n$ is a path
from $(0,0)$ to $(2n,0)$ using \emph{up} steps $(1,1)$,
\emph{down} steps $(1,-1)$, and horizontal steps of length two
$(2,0)$, and remaining weakly above the $x$-axis.

\bigskip

\begin{figure}[htb]
\begin{center}
\epsfig{file=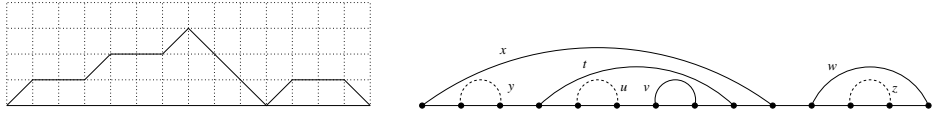, width=4.8in} \caption{A Schr\"oder path
and its representation as a two-coloured matching.}\label{seggs}
\end{center}
\end{figure}

We can represent Schr\"oder paths by using two-coloured matchings
as well. We can essentially adopt the same representation as for
Dyck paths, just using dotted arches to represent horizontal
steps. According to such a representation, dotted arches can be
contained into other arches, but they cannot contain any arch (see
figure \ref{seggs}).
This condition precisely identifies those two-coloured matchings
representing Schr\"oder paths among all two-coloured matchings.

\bigskip

Let $\mathcal{S}(n)\subseteq \mathcal{U}(n)$ denote the set of
pairs of relations $(S,R)$ on $X$ of cardinality $n$ satisfying
the following axiom:

\begin{itemize}
  \item if $xSx$, then $x$ is a minimal element for $S$.
  \hfill (\textbf{min})
\end{itemize}

Since each combinatorial interpretations of this kind of pairs of
relations is counted by Schr\"oder numbers, we call them
\emph{Schr\"oder pairs of size $n$}.

Also for axiom (\textbf{min}) an equivalent formulation can be
given which is analogous to that of proposition \ref{alt}, and
whose proof is left to the reader:
\begin{displaymath}
S\circ \mathcal{D}(S) =\mathcal{D}(S).
\end{displaymath}

Notice that, in this case, $S$ does not commute with
$\mathcal{D}(S)$; more precisely, it is
\begin{displaymath}
S\circ \mathcal{D}(S)\subseteq \mathcal{D}(S)\circ S.
\end{displaymath}

For example, referring to the matching representation of the
Schr\"oder path given in figure \ref{seggs}, we have
$y(\mathcal{D}(S)\circ S)x$, but $(y,x)\notin S\circ
\mathcal{D}(S)$.


\begin{thebibliography}{99}
\bibitem[B]{B} M. Bona,\quad \emph{Combinatorics of
permutations}\quad Discrete Mathematics and Its Applications (Boca
Raton), Chapman \& Hall, CRC, Boca Raton, FL, 2004.
\bibitem[BMCDK]{BMCDK} M. Bousquet-M\'elou, A. Claesson, M. Dukes, S. Kitaev,\quad
\emph{$(2+2)$-free posets, ascent sequences and pattern avoiding
permutations}\quad at http://arxiv.org/abs/0806.0666.
\bibitem[BEM]{BEM} A. Burstein, S. Elizalde, T. Mansour,\quad
\emph{Restricted Dumont permutations, Dyck paths, and noncrossing
partitions}\quad Discrete Math.\quad 306 (2006) 2851-2869.
\bibitem[Cl]{Cl} A. Claesson,\quad \emph{Generalized pattern avoidance}\quad
European J. Combin.\quad 22 (2001) 961-971.
\bibitem[C]{C} A. Conflitti,\quad \emph{On Whitney numbers of the order ideals of
generalized fences and crowns}\quad Discrete Math.\quad (to
appear).
\bibitem[DP]{DP} B. A. Davey, H. A. Priestley,\quad \emph{Introduction to lattices and
order}\quad Cambridge University Press, New York, 2002.
\bibitem[E]{E} S. Elizalde,\quad \emph{Multiple pattern-avoidance
with respect to fixed points and excedances}\quad Electron. J.
Combin.\quad 11 (2004) \#R51 (40pp.).
\bibitem[F]{F} T. Fine,\quad \emph{Extrapolation when very little is known about the
source}\quad Information and Control\quad 16 (1970) 331-359.
\bibitem[Fis]{Fis} P. C. Fishburn,\quad \emph{Intransitive indifference with unequal indifference intervals}
\quad J. Math. Psych.\quad 7 (1970) 144-149.
\bibitem[GP]{GP} O. Guibert, S. Pelat-Alloin,\quad \emph{Extending Fine sequences: a
link with forbidden patterns}\quad at
http://arxiv.org/abs/math/0507408v1.
\bibitem[K]{K} D. Kremer,\quad \emph{Permutations with forbidden subsequences and a
generalized Schroeder number}\quad Discrete Math.\quad 218 (2000)
121-130.
\bibitem[MaSh]{MaSh} C. L. Mallows, L. Shapiro,\quad \emph{Balls on the
lawn}\quad J. Integer Seq.\quad 2 (1999) Article 99.1.5.
\bibitem[MaSe]{MaSe} T. Mansour, S. Severini,\quad \emph{Enumeration of $(k,2)$-noncrossing
partitions}\quad Discrete Math.\quad 308 (2008) 4570-4577.
\bibitem[MM]{MM} R. J. Marsh, P. Martin,\quad \emph{Pascal arrays: counting Catalan sets}\quad
available at http://arxiv.org/abs/math.CO/0612572.
\bibitem[M]{M} J. W. Moon\quad \emph{Some enumeration problems for similarity
relations}\quad Discrete Math.\quad 26 (1979) 251-260.
\bibitem[MZ]{MZ} E. Munarini, N. Zagaglia Salvi,\quad \emph{On the rank polynomial of
the lattice of order ideals of fences and crowns}\quad Discrete
Math.\quad 259 (2002) 163-177.
\bibitem[P]{P} H. Prodinger,\quad \emph{A correspondence between ordered
trees and noncrossing partitions}\quad Discrete Math.\quad 46
(1983) 205-206.
\bibitem[STT]{STT} A. Sapounakis, I. Tasoulas, P. Tsikouras,\quad \emph{Counting strings in
Dyck paths}\quad Discrete Math.\quad 307 (2007) 2909-2924.
\bibitem[Sh]{Sh} L. W. Shapiro,\quad \emph{A Catalan
triangle}\quad Discrete Math.\quad 14 (1976), 83-90.
\bibitem[Sl]{Sl} N. J. A. Sloane,\quad \emph{The On-Line Encyclopedia of
Integer Sequences}\quad at
http://www.research.att.com/$\thicksim$njas/sequences/index.html.
\bibitem[St1]{St1} R. P. Stanley,\quad
\emph{Enumerative Combinatorics, Vol. 2}\quad Cambridge University
Press, Cambridge, 1999.
\bibitem[St2]{St2} R. P. Stanley,\quad \emph{Catalan addendum}\quad
available at
http://www-math.mit.edu/$\thicksim$rstan/ec/catadd.pdf.










\end{thebibliography}
\end{document}